\documentclass{article}
\usepackage{amsmath}
\usepackage{amscd}
\usepackage{amssymb}
\usepackage{amsthm}
\newtheorem{theorem}{Theorem}
\newtheorem{lemma}[theorem]{Lemma}
\newtheorem{corollary}[theorem]{Corollary}
\newtheorem{definition}[theorem]{Definition}

\newtheorem{proposition}[theorem]{Proposition}

\newcommand{\gl}{\lambda}
\newcommand{\gge}{\epsilon}
\newcommand{\gG}{\Gamma}
\newcommand{\ga}{\alpha}
\newcommand{\gb}{\beta}
\newcommand{\gs}{\sigma}
\newcommand{\gd}{\delta}
\newcommand{\gD}{\Delta}
\newcommand{\gal}{\text{Gal} \,}
\newcommand{\End}{\text{End} \,}
\newcommand{\tr}{\text{tr} \,}

\newcommand{\A}{\begin{pmatrix} a & b \\ c & d \end{pmatrix}}
\newcommand{\ve}{\varepsilon}

\newcommand{\mbq}{\mathbb{Q}}

\newcommand{\mbz}{\mathbb{Z}}
\newcommand{\mbf}{\mathbb{F}}
\newcommand{\mbc}{\mathbb{C}}
\newcommand{\mbn}{\mathbb{N}}
\newcommand{\mbp}{\mathbb{P}}

\newcommand{\aut}{\text{Aut}\,}
\newcommand{\ol}[1]{\overline{#1}}
\newcommand{\mc}[1]{\mathcal{#1}}

\newcommand{\mb}[1]{\mathbb{#1}}

\title{Almost all elliptic curves are Serre curves.}
\author{Nathan Jones}
\date{}
\begin{document}
\maketitle
\begin{abstract}
Using a multidimensional large sieve inequality, we obtain a bound
for the mean square error in the Chebotarev theorem for division
fields of elliptic curves that is as strong as what is implied by
the Generalized Riemann Hypothesis.  As an application we prove a
theorem to the effect that, according to height, almost all
elliptic curves are Serre curves, where a Serre curve is an
elliptic curve whose torsion subgroup, roughly speaking, has as
much Galois symmetry as possible.
\end{abstract}
\section{Introduction.} \label{intro}

Let $E$ be an elliptic curve defined over $\mbq$ and denote by
\[
\phi_{N,E} : G_\mbq \rightarrow \text{Aut}(E[N])
\]
the representation of $G_\mbq := \gal(\ol{\mbq}/\mbq)$ on the $N$-torsion of $E$.
Fixing a $\mbz/N\mbz$-basis of $E[N]$, we identify $\text{Aut}(E[N])$ with $GL_2(\mbz/N\mbz)$ and write
\[
\phi_{N,E} : G_\mbq \rightarrow GL_2(\mbz/N\mbz).
\]
The image $\phi_{N,E}(G_\mbq)$ is exactly the Galois group of the
\emph{$N$th division field of $E$ over $\mbq$}, i.e. the field
obtained by adjoining to $\mbq$ the $x$ and $y$ coordinates of the
$N$-torsion of a given Weierstrass model of $E$.
Taking the inverse limit over all $N \geq 1$ with the bases chosen compatibly,
we obtain the full torsion representation
\[
\phi_E : G_\mbq \rightarrow GL_2(\hat{\mbz}) :=
\lim_{\leftarrow}GL_2(\mbz/N\mbz).
\]
It is natural to wonder how large the image of
$\phi_E$ in $GL_2(\hat{\mbz})$ is.
\begin{definition} \label{exceptional}
The integer $N$ is said to be \emph{exceptional for
$E$} if $\phi_{N,E}$ is not surjective.
\end{definition}
To wonder about the size of the image of $\phi_E$ in $GL_2(\hat{\mbz})$ is simply
to wonder about which numbers $N$ are exceptional for $E$, and
about ``how exceptional each $N$ is,'' i.e. about the index
$[GL_2(\mbz/N\mbz):\phi_{N,E}(G_\mbq)]$.

When $E$ has complex multiplication, every $N$ except possibly
$N=2$ is exceptional, and so the image $\phi_E(G_\mbq)$ has
infinite index in $GL_2(\hat{\mbz})$. On the other hand, when $E$
does not have CM, Serre \cite{serre} has shown that the index
$[GL_2(\hat{\mbz}):\phi(G_\mbq)]$ is finite.  Equivalently, there exists an
integer $n_E$ so that
\begin{equation}  \label{n_E}
\phi_E(G_\mbq) = \pi^{-1}(\phi_{n_E,E}(G_\mbq)),
\end{equation}
where $\pi : GL_2(\hat{\mbz}) \rightarrow GL_2(\mbz/n_E\mbz)$ is
the natural projection.  In particular, this implies that any
fixed elliptic curve $E$ has only finitely many exceptional
primes, since any such exceptional prime must divide $n_E$.  One
might wonder how the integer $n_E$ (chosen minimally so that
\eqref{n_E} still holds) depends on the curve $E$.  Various results exist which bound the largest possible exceptional prime for $E$.
For example, Mazur \cite{mazur} proves that if $E$ is semistable then no prime $N \geq 11$ can be exceptional for $E$.

Define the height $H(E)$ of the elliptic curve by
\[
H(E) = \max ( |r|^3,|s|^2 ),
\]
where $r$ and $s$ are the unique integers so that $E$ has a model
of the form $y^2=x^3+rx+s$ and $\textrm{gcd}(r^3,s^2)$ is
twelfth-power free.  Duke \cite{duke} shows that, when counted
according to height, almost all elliptic curves have no
exceptional primes.  Stated precisely, he shows that if
\[
C(X) := \{ \textrm{ isomorphism classes elliptic curves } E \text{ over } \mbq : \, H(E) \leq X^6 \}
\]
and $\ve(X)$ is the set of $E \in C(X)$ which have at least one exceptional
prime, then
\begin{equation} \label{almostall}
\lim_{X \rightarrow \infty} \frac{|\ve(X)|}{|C(X)|} = 0.
\end{equation}
He does this by using a two-dimensional large sieve inequality to
prove a result which bounds the mean-square error in the
Chebotarev density theorem for the $N$th division fields of $E$
over curves of bounded height.  Using this, he shows
\begin{equation*} \label{vebound}
|\ve(X)| \ll X^4\log^BX
\end{equation*}
with an absolute (but ineffective) constant.  Since
\begin{equation} \label{sizeofCX}
C(X) = \frac{4}{\zeta(10)} X^5 + O(X^3)
\end{equation}
(c.f. \cite{brumer}), this implies \eqref{almostall}.

In \cite{grant}, Grant obtains an asymptotic formula for $\ve(X)$.  He shows that the curves which are exceptional at the
primes $2$ and $3$ contribute the main term of $|\ve(X)|$, and
that, for an explicit constant $C$,
\[
|\ve(X)| = CX^3 + O(X^{2+\gge}).
\]
for all $\gge > 0$.

This paper gives a different generalization.  The statement that an elliptic curve 
$E$ has no exceptional primes may be viewed as saying
that the Galois representation $\phi_E$ has 
``large image.'' In this paper we extend \eqref{almostall} to a result that almost all
elliptic curves have $\phi_E(G_\mbq)$ ``as large as possible.''

\section{Acknowledgment.}
This paper contains results of my Ph. D. dissertation.  I am grateful to my advisor William Duke for
his guidance.

\section{Statement of Results.}

Our main result is a theorem bounding the mean-square error in the Chebotarev theorem for division fields of 
elliptic curves.  Fix a positive integer level $N$ and a conjugacy class
\[
\mc{C} \subset GL_2(\mbz/N\mbz).
\]
We denote by
\[
\pi_E(X;N,\mc{C}) := |\{ p \leq X : \phi_{N,E}(\text{Frob}_p) \subseteq
\mc{C} \}|
\]
the function which counts the number of primes up to $X$ which are unramified in $\mbq(E[N])$ and whose Frobenius class
is contained in $\mc{C}$, and as usual
\[
\pi(X;N,d) = |\{ p \leq X : p \equiv d \mod N \}|.
\]
\begin{theorem} \label{CHEBOTAREV}
For $X \geq 1$, one has
\begin{equation*} \label{LHSCHEBOTAREV}
\frac{1}{|C(X)|}\sum_{E \in C(X)} (\pi_E(X;N,\mc{C}) -
\frac{|\mc{C}|\varphi(N)}{|GL_2(\mbz/N\mbz)|} \pi(X;N,d))^2 \ll N^8 X,
\end{equation*}
where $\varphi(N)$ denotes the Euler-phi function, and the implied constant is absolute.
\end{theorem}
In \cite{duke}, Duke proves this (without the $N^8$ factor) for prime level $N$ and the where the conjugacy
class $\mc{C}$ is replaced by a set of the form 
\[
 G_{t,d} := \{ A \in GL_2(\mbz/N\mbz) : \; \tr A = t, \; \det A = d \}
\]
Such sets are unions of conjugacy classes.  For
example, even when $N$ is prime, the set $G_{2\gl, \gl^2}$
contains two conjugacy classes, represented by the matrices
\[
\begin{pmatrix} \gl & 0 \\ 0 & \gl \end{pmatrix} \; \text{ and }
\; \begin{pmatrix} \gl & 1 \\ 0 & \gl \end{pmatrix},
\]
respectively.  Theorem \ref{CHEBOTAREV} distinguishes between
these two cases.

Our second result is an application of this theorem to the problem of counting
elliptic curves $E$ for which $\phi_E(G_\mbq)$ is as large as possible.  First of all,
how large can this image be?  Does there exist an
elliptic curve $E$ with $\phi_E$ surjective? In other words, is
there a curve $E$ with $n_E = 1$?  Serre \cite{serre} answers no.  For each
elliptic curve $E$, there is an index two subgroup $H_E \subseteq
GL_2(\hat{\mbz})$ (for a precise definition, see section
\ref{definition}) so that
\begin{equation} \label{Serre}
\phi_E(G_\mbq) \subseteq H_E.
\end{equation}
\begin{definition} \label{defofserrecurve}
We call an elliptic curve $E$ a \emph{Serre curve} when equality
holds in \eqref{Serre}.
\end{definition}
Our second theorem is
\begin{theorem} \label{main}
Let $C_{\textrm{Serre}}(X)$ denote the set
\[
\{ E \in C(X) : E \textrm{ is a Serre curve }\}.
\]
Then,
\[
\lim_{X \rightarrow \infty} \frac{|C_{\textrm{Serre}}(X)|}{|C(X)|}
= 1.
\]
\end{theorem}

In order to obtain this result that ``almost all elliptic curves
are Serre curves'', we prove an algebraic lemma which gives a
sufficient condition on an elliptic curve $E$ to be a Serre curve.
\begin{lemma} \label{serrecurvelemma}
Suppose $E$ over $\mbq$ is an elliptic
curve such that
\begin{itemize}
\item[1.]{$E$ has no exceptional primes.}

\item[2.]{$E$ is not exceptional at $4$ or $9$.}

\item[3.]{The index $[GL_2(\mbz/8\mbz) : \phi_{8,E}(G_\mbq)] \neq 2$.}

\item[4.]{There is a prime number $p > 3$ which divides the Serre
number $M_{\gD_{sf}(E)}$ (For the definition of the Serre number $M_{\gD_{sf}(E)}$, see
section \ref{definition}).}
\end{itemize}
Then, $E$ is a Serre curve.
\end{lemma}
This lemma is used together with Theorem \ref{CHEBOTAREV} to give
Theorem \ref{main}.  In a subsequent paper we plan to use Theorem \ref{main} to compute the
average value over elliptic curves of the Lang-Trotter constants, answering a question of David
and Pappalardi \cite{david}.

The paper is organized as follows:  in section \ref{boundingmeansquare} we prove Theorem \ref{CHEBOTAREV}.  Section \ref{definition}
gives the complete definition of a Serre curve, and section \ref{whichcurves} devotes itself to a proof of Lemma \ref{serrecurvelemma}.
Finally in section \ref{aaecasc} we prove Theorem \ref{main}, and in section \ref{N=4} we produce an example of a one-parameter family
of elliptic curves which are exceptional at $N=4$ but not at $N=2$.

\section{Bounding Mean-Square Chebotarev error.} \label{boundingmeansquare}

In this chapter we prove Theorem \ref{CHEBOTAREV}.  We first remark that although it 
gives a bound as strong as the appropriate Generalized Riemann Hypothesis,
the proof is unconditional.  It employs the following large sieve inequality
of Gallagher (see Lemma A of \cite{gallagher}) and follows along the same lines as
the proof of Theorem $2$ of \cite{duke}.

\begin{lemma} \label{largesieve}
For each prime number $p$ let $\Omega(p) \subseteq
(\mathbb{Z}/p\mathbb{Z})^k$ be any subset.  For each fixed $m \in
\mathbb{Z}^k$ we define
\[
P(X;m) = |\{ p \leq X : m \mod p \in \Omega(p) \} |
\]
and
\[
P(X) = \sum_{p \leq X} |\Omega(p)| p^{-k}.
\]
Let $B$ be a box in $\mathbb{R}^k$ whose sides are parallel to the
coordinate planes which has minimum width $W(B)$ and volume
$V(B)$.  If $W(B) \geq X^2$, then
\begin{equation*}
\sum_{m \in B \cap \mathbb{Z}^n}(P(X;m) - P(X))^2 \ll_k V(B)P(X).
\label{plug}
\end{equation*}
\end{lemma}
We will take $k=2$ and define the set $\Omega(p) = \Omega_\mc{C}(p)$ in such a
way that $P(X;(r,s))$ and $P(X)$ will satisfy
\begin{equation} \label{PXrs}
P(X;(r,s)) = \pi_{E_{r,s}}(X;N,\mc{C}) + O(\omega(N)). \quad
\end{equation}
and
\begin{equation} \label{PX}
P(X) = \frac{|\mc{C}|\varphi(N)}{|GL_2(\mbz/N\mbz)|} \pi(X;N,d) +
O(N^4X^{1/2}),
\end{equation}
respectively, where the implied constants are absolute.  This implies Theorem \ref{CHEBOTAREV}.

\subsection{Defining the set $\Omega_\mc{C}(p) \subseteq (\mbz/p\mbz)^2$.} \label{omega_C(p)}

We begin by quoting a result of Duke and Toth
\cite{d&t} which describes explicitly the conjugacy class in $\gal(\mbq(E[N])/\mbq)$
of the frobenius automorphism at a prime $p$ which is unramified in $\mbq(E[N])$.  The
description is given purely in terms of data attached to $E_p$, the reduction of $E$ modulo $p$.

In our context, their result may be stated as follows:  let $\mbf_p$ denote
the finite field of $p$ elements and $E_p$ any 
elliptic curve defined over $\mbf_p$.  Let $a = a(E_p) \in \mbz$ be
the trace of the Frobenius endomorphism $\phi_p$ of $E_p$ and $b = b(E_p)$ 
the index in the ring of $\mb{F}_p$-endomorphisms of $E_p$ of the 
subring generated by the frobenius endomorphism:
\[
 b = [\End_{\mb{F}_p}(E_p) : \mbz[\phi_p]].
\]

In any case (including the supersingular case), the ring $\End_{\mb{F}_p}(E_p)$ is isomorphic to an imaginary quadratic order (see Theorem $4.2$ of \cite{waterhouse}), whose discriminant we denote by $\gD = \gD(E_p)$.  The comparison of discriminants yields
\begin{equation} \label{banda}
 \gD b^2 = a^2-4p.
\end{equation}
We associate to $E_p$ the following matrix of trace $a$ and determinant
$p$:
\begin{equation} \label{sigma}
\gs(E_p) = \begin{pmatrix} (a+b\gd)/2 & b \\ b(\gD- \gd)/4 &
(a-b\gd)/2 \end{pmatrix}
\end{equation}
where for a discriminant $\gD$ we have $\gd = 0,1$ according to
whether $\gD \equiv 0,1$ mod $4$.  Because of \eqref{banda}, $\gs$ has integer entries.
\begin{theorem} \label{dtthm}
Suppose the elliptic curve $E_p$ over $(\mbz/p\mbz)$ is the reduction modulo $p$ of an elliptic curve $E$
over $\mbq$ and $p$ is prime to $N$. Then $p$ is unramified in
$\mbq(E[N])$ and the integral matrix $\gs$ defined in \eqref{sigma},
when reduced modulo $N$, represents the class of the Frobenius of
$p$ in $\textrm{Gal}\,(\mbq(E[N])/\mbq)$.
\end{theorem}
Now suppose $p > 3$ is a prime number and for $(r,s) \in
\mb{F}_p^2$, let $E_{r,s}$ denote the curve given by the equation
\begin{equation} \label{weierstrassform}
y^2 = x^3 + rx + s,
\end{equation}
and $\gD_{r,s} = -16(4r^3+27s^2)$ the
associated discriminant.  For any conjugacy class $\mathcal{C} \subset GL_2(\mbz/N\mbz)$
we define $\Omega_\mc{C}(2) = \Omega_\mc{C}(3) = \emptyset$ and for $p > 3$,
\begin{equation*} \label{Omega}
\Omega_\mathcal{C}(p) := \{ (r,s) \in \mb{F}_p^2 : \gD_{r,s} \neq
0 \text{ and } \gs(E_{r,s}) \mod N \in \mathcal{C} \}.
\end{equation*}
Since the discriminant $\gD_{r,s}$ of the curve $E_{r,s}$ is related to its minimal discriminant
$\gD$ by
\[
\gD_{r,s} = e^{12}\gD
\]
for some $e$ dividing $6$, we see from Theorem \ref{dtthm} that \eqref{PXrs} holds.  We now turn to 
verifying \eqref{PX}.

\subsection{The asymptotic in $p$ of $|\Omega_\mathcal{C}(p)|$.}

The goal of this section is to give the asymptotic of
$|\Omega_\mathcal{C}(p)|$ as $p$ ranges through the set of
prime numbers for which $\Omega_\mc{C}(p) \neq \emptyset$.  Our proof will
show that in fact,
\[
\Omega_\mc{C}(p) \neq \emptyset \Longleftrightarrow p \equiv \det \mc{C} \mod N.
\]
\begin{theorem} \label{Omegaasymptotic}
For $p$ prime congruent to $\det \mc{C}$ modulo $N$ we have
\[
|\Omega_\mathcal{C}(p)| =
\frac{|\mathcal{C}|\varphi(N)}{|GL_2(\mbz/N\mbz)|}p^2 + O(N^5p^{3/2})
\]
where the implied constant is absolute.
\end{theorem}
We observe that \eqref{PX} follows upon partial summation.  Thus, Theorem \ref{CHEBOTAREV} will follow
from Theorem \ref{Omegaasymptotic}.

To prove Theorem \ref{Omegaasymptotic}, we first express $|\Omega_\mathcal{C}(p)|$ in terms of a weighted class number.
Define the set
\[
\mc{T}_\mc{C}(p) := \{ A \in M_{2 \times 2}(\mbz) : \; \det A = p, \, A \mod N \in \mc{C} \},
\]
and the subset of elliptic matrices
\[
\mc{T}^e_\mc{C}(p) := \{ A \in \mc{T}_\mc{C}(p) : \; ( \tr A )^2 - 4 \det A < 0 \}.
\]
Since $\mc{C}$ is stable by $SL_2(\mbz/N\mbz)$-conjugation, both of the above sets are stable by $SL_2(\mbz)$-conjugation.

Note:  Throughout the rest of this paper we will use the standard notation
\[
 \gG(1) := SL_2(\mbz).
\]
\begin{proposition} \label{classnumberprop}
\[
|\Omega_\mc{C}(p)| = \frac{p-1}{2} \sum_{\ga \in \mc{T}_\mc{C}^e(p) \,//\, \gG(1)} \frac{1}{|\gG(1)_\ga|},
\]
where $\mc{T}_\mc{C}^e(p) \,//\, \gG(1)$ is the set of $\gG(1)$-conjugation orbits in $\mc{T}_\mc{C}^e(p)$ and
\[
\gG(1)_\ga := \{ \gamma \in \gG(1) : \; \gamma \ga = \ga \gamma \}.
\]
\end{proposition}
This proposition, together with the following lemma, imply Theorem \ref{Omegaasymptotic}.
\begin{lemma} \label{traceformulalemma}
If $p \equiv \det \mc{C} \mod N$ then
\[
\sum_{\ga \in \mc{T}^e_\mc{C}(p) \,//\, \gG(1)} \frac{1}{|\gG(1)_\ga|} = \frac{2 |\mc{C}|}{|SL_2(\mbz/N\mbz)|} p + O(N^5p^{1/2}),
\]
with an absolute constant.
\end{lemma}
\begin{proof}
This is Corollary 8 of \cite{jones}
\end{proof}
The remainder of this section is devoted to proving Proposition \ref{classnumberprop}.  We note that 
\[
\Omega_\mc{C}(p) = \{ (r,s) \in (\mbz/p\mbz)^2 : \; \gD_{r,s} \neq 0 \; \text{ and } \; \gs(E_{r,s}) \in \mc{T}^e_\mc{C}(p) \}.
\]
At this point we must give a finer description of $\mc{C}$.
For any divisor $M$ of $N$ and integers $\ol{T}$, $\ol{D}$ modulo $N/M$, define
\[
\mc{T}_{N/M}(\ol{T},\ol{D}) = \{ A \in M_{2 \times 2}(\mbz/(N/M)\mbz) : \; (\tr A, \det A) \equiv (\ol{T},\ol{D}) \mod N/M \}
\]
and
\[
\mc{T}^*_{N/M}(\ol{T},\ol{D}) = \{ A \in \mc{T}_{N/M}(\ol{T},\ol{D}) : \; A \text{ is non-scalar mod each prime } l \mid N/M \}.
\]
The following lemma is a corollary of Proposition $7$ of \cite{jones} describing the structure of conjugacy classes in the group
$GL_2(\mbz/N\mbz)$.
\begin{lemma} \label{description}
Any conjugacy class 
\[
\mc{C} \subset GL_2(\mbz/N\mbz)
\]
has the form
\begin{equation*} \label{conjugacy}
\mc{C} = \gl I + M \mc{T}^*_{N/M}(\ol{T},\ol{D}),
\end{equation*}
where $\gl$ is an integer satisfying $0 \leq \gl < M$.
\end{lemma}

We would like to
partition $\mc{T}^e_\mc{C}(p)$ into subsets which are stable by $\gG(1)$-conjugation.  Let $\mc{T}^*(T,D,f)$ denote
\[
\{ A = \A \in M_{2 \times 2}(\mbz) : \; \tr A = T, \, \det A = D, \, \gcd(b,d-a,c) = f \}.
\]
We note then that the trace $t$ and determinant $d$ of any matrix in the set $\gl I + M \mc{T}^*(T,D,f)$ satisfy
\begin{equation} \label{traceofA}
t = 2\gl + MT, \quad d = \gl^2 + M \gl T + M^2 D, \quad \text{ and } \quad t^2-4d = M^2 \left( T^2 - 4D \right).
\end{equation}
Thus, from Lemma \ref{description} we see that
\[
\mc{T}^e_\mc{C}(p) = \bigsqcup_{(T,D)}\bigsqcup_{{\begin{substack} {f \geq 1 \\ \gcd(f,N/M) = 1} \end{substack}}} \left( \gl I + M\mc{T}^*(T,D,f) \right),
\]
where $(T,D)$ runs over integer pairs satisfying
\begin{equation*} \label{TandD}
(T,D) \equiv (\ol{T},\ol{D}) \mod N/M, \; p = \gl^2 + M\gl T + M^2 D,\, \text{ and } \, (2\gl + MT)^2 < 4p.
\end{equation*}
Defining $\Omega^*(\gl,M,T,D,f)$ by
\[
\{ (r,s) \in (\mbz/p\mbz)^2 : \; \gD_{r,s} \neq 0 \text{ and } \gs(E_{r,s}) \in \gl I + M\mc{T}^*(T,D,f) \},
\]
Proposition \ref{classnumberprop} is reduced to showing that
\begin{equation} \label{newsum}
|\Omega^*(\gl,M,T,D,f)| = \frac{p-1}{2} \sum_{\ga \in \left( \gl I + M\mc{T}^*(T,D,f) \right) \,//\, \gG(1)} \frac{1}{|\gG(1)_\ga|}.
\end{equation}

\begin{lemma}
$\Omega^*(\gl,M,T,D,f)$ is equal to 
\[
\{ (r,s) \in (\mbz/p\mbz)^2 : \; \gD_{r,s} \neq 0, \; b(E_{r,s}) = Mf \; \text{ and } \; a(E_{r,s}) = 2\gl + MT \}.
\]
\end{lemma}
\begin{proof}
The containment ``$\Omega^*(\gl,M,T,D,f) \subseteq$ \dots'' is immediate from \eqref{sigma} and \eqref{traceofA}. 
The reverse containment comes from the fact that, for fixed $t$ and $p$, the two equations
\[
t = 2\gl + MT \quad \text{ and } \quad p = \gl^2 + M \gl T + M^2 D
\]
have a unique solution $(\gl,T,D) \in \{ 0,1, \dots,M-1\} \times \mbz^2$, if they have one at all.  This fact is
immediate when $M$ is odd.  If $M$ is even, we see from the first equation that the only way two distinct solutions can exist is if one solution looks like $(\gl,T,D)$ with $\gl \in \{0,1,\dots,M/2-1\}$ and the other solution has the form $(\gl+M/2,T-1,D')$ for some integer $D'$.  But then the second equation gives us the contradiction that
\[
\gl^2 + M \gl T - p \equiv 0 \mod M^2 \quad \text{ and } \quad \gl^2 + M \gl T - p \equiv \frac{M^2}{4} \left( 1-2T \right) \mod M^2.
\]
\end{proof}

We now summarize some fundamental facts about imaginary quadratic orders.  More details may be found, for example, in \S $7$ of \cite{cox}.  An \emph{imaginary quadratic order} $\mc{O}$ is a subring (containing $1$) of an imaginary quadratic field $K$ which contains a basis of $K$ over $\mbq$ and has rank $2$ as an free abelian group.  For each negative number $\gD$ satisfying
\[
 \gD \equiv 0 \text{ or } 1 \mod 4,
\]
there is a unique imaginary quadratic order of discriminant $\gD$, which we will denote by $\mc{O}(\gD)$.  Orders $\mc{O}(\gD')$ which contain a given order $\mc{O}(\gD)$ are exactly those orders whose discriminant $\gD'$ satisfies
\[
 f^2 \gD' = \gD, \quad\quad f = [\mc{O}(\gD') : \mc{O}(\gD)].
\]
Every imaginary quadratic order $\mc{O}$ is contained in a unique maximal imaginary quadratic order
\[
 \mc{O} \subseteq \mc{O}_{\max} = \mc{O}_K \subset K,
\]
which is the ring of integers of $K$.  The ideal class group $\mc{C}(\mc{O})$ is the group of invertible fractional ideals of $\mc{O}$ modulo the subgroup of principal fractional ideals.  This is a finite group whose size we denote by $h(\mc{O})$.
\begin{lemma} \label{ellipticcurveswithO}
Suppose $p \geq 5$ is prime and $t$ is any integer satisfying $t^2<4p$.  Let $\mc{O}$ be any imaginary quadratic order containing the order of discriminant $t^2-4p$.  The number of elliptic curves $E_{r,s}$ over $\mb{F}_p$ of the form \eqref{weierstrassform} which satisfy 
\[
\tr (\phi_p) = t \quad \text{ and } \quad \End_{\mbf_p}(E_{r,s}) = \mc{O}
\]
is given by
\[
 \frac{p-1}{|\mc{O}^*|} h(\mc{O}),
\]
where $\mc{O}^*$ is the group of units of $\mc{O}$.
\end{lemma}
\begin{proof}
The following theorem restates Theorems $4.2$ and $4.5$ of \cite{waterhouse}, specialized to our situation.  See also \cite{schoof}, which corrects a small error in the original proof.  The original work is due to Deuring \cite{deuring}.
\begin{theorem}
Let $t$ be any integer satisfying $ t^2 < 4p$.  Then the following are precisely the rings which occur as rings of $\mb{F}_p$-endomorphisms of some elliptic curve $E_p$ over $\mb{F}_p$ satisfying $a(E_p) = t$:
\begin{itemize}
 \item if $t \neq 0$:  all complex quadratic orders containing $\mc{O}(t^2-4p)$;
 \item if $t=0$: all complex quadratic orders $\mc{O}$ satisfying
\[
 \mc{O}(-4p) \subset \mc{O} \quad \text{ and } \quad p \nmid [\mc{O}_{\max} : \mc{O}].
\]
\end{itemize}
Furthermore, given such an order $\mc{O}$, the number of $\mb{F}_p$-isomorphism classes of elliptic curves $E_p$ over $\mb{F}_p$ satisfying 
\[
 a(E_p) = t \quad \text{ and } \quad \End_{\mb{F}_p}(E_p) = \mc{O}
\]
is equal to $h(\mc{O})$.
\end{theorem}
Note that, since $p \geq 5$, every $\mbf_p$-isomorphism class contains an elliptic curve of the form \eqref{weierstrassform}.  By the theorem, the proof of Lemma \ref{ellipticcurveswithO} is reduced to showing that whenever $E_{r,s}$ is the form \eqref{weierstrassform} with $\tr \phi_p = t$ and $\End_{\mbf_p}(E_{r,s}) = \mc{O}$, the number of elliptic curves of the same form which are isomorphic over $\mbf_p$ to $E_{r,s}$ is $(p-1)/|\mc{O}^*|$.
Such elliptic curves are exactly those given by the equations
\[
E_{ru^4,su^6} : \; y^2 = x^3 + ru^4x + su^6, \quad u \in (\mbz/p\mbz)^*.
\]
In case $|\mc{O}^*| = 2$, the $j$-invariant $j(E_{r,s})$ cannot be equal to $0$ or $1728$, i.e. neither $r$ nor $s$ can be equal to zero.  In this case, $E_{ru^4,su^6} = E_{r(u')^4,s(u')^6}$ if and only if $u = \pm u'$ and we count exactly $(p-1)/2$ distinct $E_{ru^4,su^6}$'s.  The case $|\mc{O}^*| = 4$ occurs exactly when $\mc{O} = \mc{O}(-4) = \mbz[i]$ is the ring of Gaussian integers, and this happens only if $j(E_{r,s}) = 1728$ and $s=0$.  Since
\[
 \mc{O}(t^2-4p) \subset \mc{O}(-4),
\]
we see by relating the discriminants that $t$ must be even and that $p \equiv 1 \mod 4$.  Choosing $i_p \in (\mbz/p\mbz)^*$ satisfying $i_p^2 = -1$, we note that in this case $E_{ru^4,su^6} = E_{r(u')^4,s(u')^6}$ if and only if $u/u' \in \{\pm i_p, \pm 1 \}$, and so there are again exactly $(p-1)/|\mc{O}^*|$ elliptic curves of the form \eqref{weierstrassform} isomorphic over $\mbf_p$ to $E_{r,s}$.  The $j(E_{r,s}) = 0$ case is quite similar, so we omit it.  This finishes the proof of Lemma \ref{ellipticcurveswithO}.
\end{proof}
Returning to the verification of \eqref{newsum}, we see by the two lemmas and \eqref{traceofA} that
\[
|\Omega^*(\gl,M,T,D,f)| = \frac{p-1}{|\mc{O}\left(\frac{T^2-4D}{f^2}\right)^*|} h\left( \mc{O} \left( \frac{T^2-4D}{f^2} \right) \right).
\]

Now we use a theorem which equates the counting of weighted $\gG(1)$-orbits of matrices of a fixed trace and determinant (of negative discriminant) with the counting of weighted ideal classes in the imaginary quadratic order of the same discriminant.  We denote by $Q^*(\gD)$ the set of primitive integral binary quadratic forms of discriminant $\gD$ and $Q^*_+(\gD)$ the subset of positive definite forms, both acted on by the classical $\gG(1)$-action
\[
f \cdot \begin{pmatrix} a & b \\ c & d \end{pmatrix} (x,y) = f(ax+by,cx+dy).
\]
By $Q^*(\gD) \,//\,\gG(1)$ and $Q^*_+(\gD) \,//\, \gG(1)$ we denote the corresponding orbit spaces under this action.
\begin{theorem}
Let $T$ and $D$ be integers and $f$ a positive integer satisfying
\[
T^2-4D < 0 \quad \text{ and } \quad \frac{T^2-4D}{f^2} \in \mbz, \quad \frac{T^2-4D}{f^2} \equiv 0,1 \mod 4.
\]
Then there are set bijections
\[
\mc{T}^*(T,D,f)\,//\,\gG(1) \longleftrightarrow Q^*\left( \frac{T^2-4D}{f^2} \right) \,//\,\gG(1)
\]
and
\[
Q^*_+\left( \frac{T^2-4D}{f^2} \right) \,//\,\gG(1) \longleftrightarrow \mc{C}\left(\mc{O} \left(\frac{T^2-4D}{f^2} \right) \right),
\]
\end{theorem}
\begin{proof}
We first observe that, whenever $\mc{T}^*(T,D,f) \neq \emptyset$ (which is equivalent to the three given conditions), there are unique integers $T'$, $D'$ and $\gl \in \{ 0,1, \dots, f-1 \}$ so that
\[
\mc{T}^*(T,D,f) = \gl I + f \mc{T}^*(T',D',1).
\]
Since $T^2-4D = f^2((T')^2-4D')$, the first bijection in the theorem is induced by the bijection
\[
\mc{T}^*(T',D',1) \longleftrightarrow Q^*\left( (T')^2-4D' \right)
\]
given by sending the matrix $\begin{pmatrix} a & b \\ c & d \end{pmatrix}$ to the form $cx^2 + (d-a)xy -by^2$ and the form $\ga x^2 + \gb xy + \gamma y^2$ to the matrix $\begin{pmatrix} (t-\gb)/2 & -\gamma \\ \ga & (t-\gb)/2 \end{pmatrix}$.
The second bijection is Theorem 7.7 in \cite{cox}.
\end{proof}

We observe that for any matrix $\ga \in \mc{T}^*(T,D,f)$, we have
\[
|\gG(1)_\ga| = |\mc{O}\left(\frac{T^2-4D}{f^2} \right)^*|,
\]
and the common value can be greater than $2$ only when $\frac{T^2-4D}{f^2} \in \{-3,-4\}$, in which case $h(\mc{O}\left(\frac{T^2-4D}{f^2} \right)) = 1$.  We conclude:

\begin{corollary}
\[
\frac{2}{|\mc{O}\left( \frac{T^2-4D}{f^2} \right)^*|} h\left( \mc{O}\left( \frac{T^2-4D}{f^2} \right) \right) = \sum_{\ga \in \left( \gl I + M\mc{T}^*(T,D,f) \right)\,//\,\gG(1) } \frac{1}{|\gG(1)_\ga|}.
\]
\end{corollary}
By the corollary, \eqref{newsum} follows and we have proved Proposition \ref{classnumberprop}.

\section{The definition of a Serre curve.} \label{definition}

We now describe the subgroup $H_E$ mentioned in Definition \ref{defofserrecurve},
following the discussion proceeding Proposition $22$ of \cite{serre}.
Suppose that $E$ is given by the equation
\[
 y^2 = x^3 + rx + s = (x-e_1)(x-e_2)(x-e_3).
\]
Then $\{e_1,e_2,e_3\}$ is the set of $x$-coordinates of
the non-trivial 2-torsion of $E$.  The discriminant $\Delta$ of this model of $E$ is given by
\begin{equation} \label{Delta}
\Delta = ((e_1-e_2)(e_1-e_3)(e_2-e_3))^2.
\end{equation}
Thus, one has
\begin{equation*} \label{contain2}
\mbq(\sqrt{\gD}) \subset \mbq(E[2]).
\end{equation*}
Because of the action of $\aut{E[2]} \simeq GL_2(\mbz/2\mbz)$ on the $e_i$'s
we have a group isomorphism between $GL_2(\mbz/2\mbz)$ and the symmetric group on 
three letters:
\[
GL_2(\mbz/2\mbz) \simeq S_3.
\]
By \eqref{Delta} we see that for any Galois automorphism $\gs \in \gal(\mbq(E[2])/\mbq) \subset S_3$, 
\begin{equation} \label{epsilon}
\gs : \sqrt{\gD} \mapsto \ve(\gs) \sqrt{\gD},
\end{equation}
where $\ve$ denotes the signature character on $S_3$.  In particular we note that if $\sqrt{\gD} \in \mbq$ then
\[
\gal(\mbq(E[2])/\mbq) \subset A_3 = \text{ the alternating group on $3$ letters.}
\]
In this case we define the \emph{Serre number} $M_{1}$ to be $2$ and the \emph{Serre subgroup} $H_2$ by
\[
 H_2 := A_3.
\]
Otherwise, $\mbq(\sqrt{\gD})$ is a quadratic extension, which in particular is abelian.  
Since each abelian extension of $\mathbb{Q}$ is contained in a
cyclotomic extension, one may choose a positive integer $D$ so that
\begin{equation*} \label{contain}
\mbq(\sqrt{\Delta}) \subset \mbq(\zeta_{D}) \subset \mbq(E[D]),
\end{equation*}
where as usual $\zeta_D$ denotes a primitive $D$-th root of unity and the 
second containment comes from the Weil pairing (c.f. \cite{silverman}, for example).
\begin{lemma} \label{smallestD}
Let $W$ be any square-free integer and define the positive integer
$D_W$ by
\[
D_W = \begin{cases} |W| & \text{ if } W \equiv 1 \mod 4 \\
                    4|W| & \text{ otherwise.}
\end{cases}
\]
Then we have
\begin{equation*} 
\mbq(\sqrt{W}) \subset \mbq(\zeta_D) \Leftrightarrow D_W \textrm{ divides } D.
\end{equation*}
Furthermore, for such a $D$ and $\gs \in \gal(\mbq(E[D])/\mbq) \subseteq GL_2(\mbz/D\mbz)$, we have
\begin{equation} \label{Kronecker}
\gs : \sqrt{W} \mapsto \left( \frac{W}{\det \gs} \right) \sqrt{W}.
\end{equation}
Here we use the Kronecker symbol $\begin{displaystyle} \left( \frac{W}{\cdot} \right) \end{displaystyle} := \left( \frac{W/|W|}{\cdot} \right) \cdot \prod_{p \mid W} \left( \frac{p}{\cdot} \right)$, where
\[
\left( \frac{2}{\cdot} \right) := (-1)^{((\cdot)^2 - 1)/8}, \quad \text{ and } \quad \left( \frac{\pm 1}{\cdot} \right) = \left( \pm 1 \right)^{((\cdot)-1)/2}.
\]
\end{lemma}
\begin{proof}
These are standard results from algebraic number theory, together with Theorem $6.6$ of \cite{shimura}.
\end{proof}
By the lemma we see that
\begin{equation} \label{conditiononD}
\mbq(\sqrt{\Delta}) \subset \mbq(\zeta_{D}) \Leftrightarrow D_{\gD_{sf}}
\textrm{ divides } D,
\end{equation}
where $\gD_{sf} = \gD_{sf}(E)$ is the square-free part of the discriminant $\gD$
of $E$.  For any square-free number $W$ we define the ``Serre number'' 
\[
M_{W} = \begin{cases} 2|W| & \text{ if } W \equiv 1 \mod 4 \\
4|W| & \text{ otherwise } \end{cases},
\]
to be the least common multiple of $2$ and $D_{W}$.  Thus in particular, $\mbq(E[M_{\gD_{sf}}])$ is the compositum of $\mbq(E[2])$ and $\mbq(E[D_{\gD_{sf}}])$.  We furthermore define the subgroup $H_{M_W}$ by
\[
 H_{M_W} = \ker \left(\left( \frac{W}{\det(\cdot)} \right) \ve(\cdot)\right) \subset GL_2(\mbz/M_W\mbz),
\]
where here we have extended the definition of the signature character $\ve$ in the natural way to any even level:
\begin{equation} \label{ve}
 \ve : GL_2(\mbz/2m\mbz) \longrightarrow GL_2(\mbz/2\mbz) \longrightarrow \{ \pm 1 \}.
\end{equation}
Later in the paper we will casually refer to ``$\ker \ve$'', hoping that in each instance
its domain will be clear from context.

By virtue of \eqref{epsilon} and \eqref{Kronecker}, we see that
\[
 \gal(\mbq(E[M_{\gD_{sf}}]/\mbq) \subseteq H_{M_{\gD_{sf}}}.
\]
The subgroup $H_E$ of $GL_2(\hat{\mbz})$ referred to in \eqref{Serre} is simply
\[
H_E = \pi_{M_{\gD_{sf}}}^{-1}(H_{M_{\gD_{sf}}}),
\]
where $\pi_{M_{\gD_{sf}}} : GL_2(\hat{\mbz}) \longrightarrow GL_2(\mbz/M_{\gD_{sf}}\mbz)$ is the natural projection.
$H_E$ is evidently an index 2 subgroup of $GL_2(\hat{\mbz})$ and
\[
\phi_E(G_\mbq) \subseteq H_E.
\]
An elliptic curve $E$ is a Serre curve if
$\phi_E(G_\mbq) = H_E$. In other words,
an elliptic curve is a Serre curve exactly when, for every integer
$m$, we have
\[
[GL_2(\mbz/m\mbz) : \phi_{m,E}(G_\mbq) ] =
\begin{cases}   2 & \text{ if } M_{\gD_{sf}(E)} \mid m \\
                1 & \text{ otherwise.}
\end{cases}
\]
We will refer to $H_{M_{\gD_{sf}(E)}} \subset GL_2(\mbz/M_{\gD_{sf}(E)}\mbz)$ (and by abuse of notation,
also to $H_E \subset GL_2(\hat{\mbz})$) as the ``Serre subgroup
associated to $E$.''

\section{Which elliptic curves are Serre curves?} \label{whichcurves}

If $N$ is exceptional for $E$ (see Definition \ref{exceptional}) then so is any
multiple of $N$. 
\begin{definition} \label{minimalexceptional}
 We call an integer $N$ \emph{minimal exceptional
for $E$} if it is exceptional for $E$ and none of its proper nontrivial
divisors are exceptional for $E$.
\end{definition}
For example, if $E$ is a Serre
curve, then the Serre number $M_{\gD_{sf}(E)}$ (see section \ref{definition}) is a
minimal exceptional number for $E$.  Also, any exceptional prime
$p$ of $E$ is minimal exceptional.

The proof of Lemma \ref{serrecurvelemma} uses only the theory
of the groups $GL_2(\mbz/N\mbz)$ (especially for $N$ divisible by $2$ and $3$, complimenting \cite{serre2})
as well as a few facts about cyclotomic fields.  The arguments are similar to those given in
Kani's appendix to \cite{cojocaru2}.  Two separate issues arise:
(1) which numbers $N$ can actually occur as minimal exceptional
numbers for an elliptic curve and (2) the stability of the Serre
number $M_{\gD_{sf}(E)}$.  We treat them in that order.

We will make repeated use of
\begin{lemma} \label{commutator}
The commutator subgroup $(GL_2(\mbz/p^n\mbz))'$ of
$GL_2(\mbz/p^n\mbz)$ is given by
\begin{equation*}
(GL_2(\mbz/p^n\mbz))' = \begin{cases}   SL_2(\mbz/p^n\mbz) & \text{ if } p \neq 2 \\
\ker(\ve) \cap SL_2(\mbz/2^n\mbz) & \text{
if } p = 2 \end{cases} \label{comm}
\end{equation*}
(see \eqref{ve}.)  For $p \geq 5$, the group $SL_2(\mbz/p^n\mbz)$ is equal to its own commutator:
\[
 (SL_2(\mbz/p^n\mbz))' = SL_2(\mbz/p^n\mbz) \quad \quad ( p \geq 5 ).
\]
\end{lemma}

\subsection{Minimal exceptional numbers of elliptic
curves.} \label{minimalexception}

The following lemma gives us a restriction on which positive
integers $N$ can occur as a minimal exceptional number of an
elliptic curve.  Throughout the remainder of the paper we will 
sometimes use the abbreviation
\[
G_N := \gal(\mbq(E[N])/\mbq),
\]
suppressing the dependence on the elliptic curve $E$.
\begin{lemma} \label{minimal}
Let $E$ be an elliptic curve over $\mbq$.  Suppose that $N \in
\mbn$ is minimal exceptional for $E$.  Then,
\[
N \in \{\text{ prime numbers } \} \cup \{ M_{\gD_{sf}(E)} \} \cup \{4, 8, 9\}.
\]
If $8$ is  a minimal exceptional number for $E$, then
\[
 [GL_2(\mbz/8\mbz) : \phi_{E,8}(G_\mbq)] = 2.
\]
\end{lemma}
\begin{proof}
Let us assume that $N$ is not prime.  If $N$ is exceptional for
$E$ then we have
\[
G_N \subsetneq GL_2(\mbz/N\mbz).
\]
If $N$ is \emph{minimal} exceptional, we have $G_d =
GL_2(\mbz/d\mbz)$ for each proper divisor $d$ of $N$.  Therefore
the canonical map
\begin{equation}
G_N \twoheadrightarrow GL_2(\mbz/d\mbz), \quad \forall d \mid N
\label{ontod}
\end{equation}
is a surjection. By the surjectivity of the  Weil pairing, we also
see that the determinant map
\begin{equation}
\det : G_N \twoheadrightarrow (\mbz/N\mbz)^* \label{ontodet}
\end{equation}
is surjective.  We consider the question: for which composite
numbers $N$ does there exist a subgroup $G_N$ of
$GL_2(\mbz/N\mbz)$ satisfying conditions \eqref{ontod} and
\eqref{ontodet}?  We divide the investigation into cases according
to whether $N$ is a prime power or not.  We tackle the latter case
first:

\textbf{Case 1}:  $N$ is not a prime power.  Let $p$ be the
smallest prime divisor of $N$.  Suppose $p^n || N$ and write $M :=
N/p^n \, (\neq 1)$.  By Galois theory we must have
\[
\mbq \subsetneq \mbq(E[p^n]) \cap \mbq(E[M]) =: F
\]
Let $H = \gal(F/\mbq)$.  If $H$ is not simple, replace it by any nontrivial simple
quotient, and replace $F$ by the corresponding field. The next lemma is a corollary of the discussion on page IV-$25$ of \cite{serre2}.
\begin{lemma} \label{jordanhoelder}
 If $N_1$ and $N_2$ are relatively prime positive integers than the groups $GL_2(\mbz/N_1\mbz)$ and $GL_2(\mbz/N_2\mbz)$
have no common simple nonabelian quotient.
\end{lemma}
Since $H$ is a common quotient of the groups
\[
 \gal(\mbq(E[M])/\mbq) = GL_2(\mbz/M\mbz) \quad \text{ and } \quad \gal(\mbq(E[p^n])/\mbq) = GL_2(\mbz/p^n\mbz)
\]
we conclude that $H$ is abelian.  From this and Lemma \ref{commutator} it follows that
\[
F \subset \mbq (\zeta_{M}).
\]
If $p > 2$ then we must similarly have $F \subset \mbq(\zeta_{p^n})$.  Since
\begin{equation} \label{disjointness}
 \mbq(\zeta_M) \cap \mbq(\zeta_{p^n}) = \mbq,
\end{equation}
we conclude that $F = \mbq$, contradicting that $H$ is nontrivial.  Therefore
we must have $p = 2$.  But then using Lemma \ref{commutator} we
similarly conclude that
\[
\mbq \neq F \subset \mbq (\sqrt{\Delta_E},\zeta_{2^n}) \cap \mbq
(\zeta_{M}).
\]
If $n\leq 1$ then we must have $F = \mbq(\sqrt{\gD_E})$, and we
see that $N$ is a multiple of the Serre number $M_{\gD_{sf}(E)}$.  If $n \geq
2$ then we reason as follows: since the Galois group
$\textrm{Gal}(\mbq (\sqrt{\Delta_E},\zeta_{2^n})/\mbq)$ has order
a power of two, $F$ must contain a quadratic subfield
$K$.  By \eqref{disjointness}, we conclude
that if $n=2$, $K$ must be one of the fields
\[
\mbq(\sqrt{\gD_E}), \, \mbq(\sqrt{-\gD_E})
\]
and if $n \geq 3$ that $K$ must be one of the fields
\[
\mbq(\sqrt{\gD_E}), \, \mbq(\sqrt{-\gD_E}), \,\mbq(\sqrt{2\gD_E}),
\, \mbq(\sqrt{-2\gD_E}).
\]
Thus in any case by \eqref{conditiononD}, $N$ is a multiple of the
Serre number of $E$, which implies that $N$ \emph{is} the Serre
number of $E$, since $N$ is assumed to be minimal exceptional.  We
have shown that the Serre number of $E$ is the only minimal
exceptional number which is not a prime power.

\textbf{Case 2}:  $N = p^n$ is a prime power with $n \geq 2$.  If $p \geq 5$, we reason as follows:
Taking commutators of \eqref{ontod} we have a surjection
\[
 (G_{p^n}(E))' \twoheadrightarrow SL_2(\mbz/p^{n-1}\mbz) = (GL_2(\mbz/p^{n-1}\mbz))'.
\]
By Lemma $3$ on page IV-23 of \cite{serre2}, this implies that 
$(G_{p^n}(E))' = SL_2(\mbz/p^n\mbz)$.  But now since 
\[
SL_2(\mbz/p^n\mbz) \subset G_{p^n}(E)
\]
we conclude by \eqref{ontodet} that $G_{p^n}(E) = GL_2(\mbz/p^n\mbz)$, contradicting the fact that $p^n$ is
exceptional.  We conclude that $p \in \{2,3\}$.

Now consider the exact sequence
\[
1 \rightarrow K \rightarrow GL_2(\mbz/p^n\mbz) \rightarrow
GL_2(\mbz/p^{n-1}\mbz) \rightarrow 1.
\]
Here $K = I + p^{n-1}M_2(\mbz/p\mbz)$.  Since $G_{p^n}$ surjects
onto $GL_2(\mbz/p^{n-1}\mbz)$, we have the exact sequence
\[
1 \rightarrow K\cap G_{p^n} \rightarrow G_{p^n} \rightarrow
GL_2(\mbz/p^{n-1}\mbz) \rightarrow 1.
\]
First we show that if $n \geq 3$ then
\begin{equation}
I + p^{n-1}\{ \text{ traceless matrices } \} \subseteq K\cap
G_{p^n}. \label{traceless}
\end{equation}
This is seen by choosing any preimage
\[
\begin{pmatrix} 1 & p^{n-2} \\ 0 & 1 \end{pmatrix} + p^{n-1}A \in G_{p^n}
\]
of the matrix $\begin{pmatrix} 1 & p^{n-2} \\ 0 & 1 \end{pmatrix}
\in GL_2(\mbz/p^{n-1}\mbz)$ and observing that, if $n\geq 3$,
\begin{equation} \label{congruence}
\left(
\begin{pmatrix} 1 & p^{n-2} \\ 0 & 1 \end{pmatrix} + p^{n-1}A \right)
^p \equiv \begin{pmatrix} 1 & p^{n-1} \\ 0 & 1
\end{pmatrix} \mod p^n,
\end{equation}
which shows that the matrix $I + p^{n-1}\begin{pmatrix} 0 & 1 \\ 0
& 0 \end{pmatrix} \in K \cap G_{p^n}$.  Now let $\begin{pmatrix} a
& b \\ c & d \end{pmatrix}$ be any matrix in $GL_2(\mbz/p\mbz)$
and choose a matrix $A \in G_{p^n}(E)$ with
\[
A \equiv \begin{pmatrix} a & b \\ c & d \end{pmatrix} \mod p.
\]
We then have
\[
A(I + p^{n-1}\begin{pmatrix} 0 & 1 \\ 0 & 0 \end{pmatrix})A^{-1} =
I + \frac{1}{ad-bc}p^{n-1}\begin{pmatrix} -ac & a^2 \\ -c^2 & ac
\end{pmatrix} \in K \cap G_{p^n}(E).
\]
Letting the matrix $\begin{pmatrix} a & b \\ c & d \end{pmatrix}$
vary modulo $p$ we see that \eqref{traceless} holds.  From this we
see that $G_{p^n}$ must be an index 2 subgroup of
$GL_2(\mbz/p^n\mbz)$.  Thus, there is a character
\begin{equation} \label{character}
\chi : GL_2(\mbz/p^n\mbz) \rightarrow \{ \pm 1 \} \quad \text{ with } \quad G_{p^n}(E) = \ker \chi.
\end{equation}
By Lemma \ref{commutator}, we
see that if $p=3$ then $\chi$ factors through the determinant,
i.e. that there is a character
\[
\gd : (\mbz/3^n\mbz)^* \rightarrow \{\pm 1\}
\]
with $\chi = \gd \circ \det$.  This implies that
$SL_2(\mbz/3^n\mbz) \subset G_{3^n}(E)$, which says that
$G_{3^n}(E) = GL_2(\mbz/3^n\mbz)$ by \eqref{ontodet}, a
contradiction.  Thus, the only (non-prime) power of $3$ which can
be minimal exceptional for $E$ is $9$.

Now let us return to \eqref{character} with $p=2$.  In this case
Lemma \ref{commutator} says that either $\chi$ or $\chi \ve$
factors through the determinant, according to whether
\[
\chi( \begin{pmatrix} 1 & 1 \\ 0 & 1
\end{pmatrix} ) = 1 \text{ or } -1.
\]
Now if $\chi(\begin{pmatrix} 1 & 1 \\ 0 & 1 \end{pmatrix}) = 1$
then $SL_2(\mbz/2^n\mbz) \subseteq G_{2^n}$, a contradiction.  Thus we must
have $\chi(\begin{pmatrix} 1 & 1 \\ 0 & 1 \end{pmatrix}) = -1$.
Therefore since $\chi \ve$ factors through the determinant,
we have
\[
\chi = \ve \cdot (\gd \circ \text{det}),
\]
where $\gd : (\mbz/2^n\mbz)^* \rightarrow \{ \pm 1 \} $ is a
character.  Now pick $X = \begin{pmatrix} 1 & 1 \\ 0 & 1
\end{pmatrix} + 2^{n-1}A \in G_{2^n}$.  We have $\det X = 1$
or $1+2^{n-1}$.  One verifies by induction that for $n \geq 3$,
\[
1 + 2^{n-1} \equiv 5^{2^{n-3}} \mod 2^n,
\]
so for $n > 3$ we must have $\gd(\det X) = 1$, contradicting \eqref{character}. We
have shown that for any elliptic curve $E$, if $N \neq M_{\gD_{sf}(E)}$ is a
composite minimal exceptional number for $E$, then $N = 4$, $8$ or
$9$, where if $N=8$ there is a real character
\[
\chi : GL_2(\mbz/8\mbz) \rightarrow \{\pm 1\}
\]
with $G_8(E) = \ker \chi$. This concludes the proof of
Lemma \ref{minimal}.
\end{proof}

\subsection{Stability of the Serre number $M_{\gD_{sf}(E)}$.}

Continuing the proof of Lemma \ref{serrecurvelemma}, we will now show
that under the assumptions stated therein, we have
\begin{equation} \label{HME}
G_{M_{\gD_{sf}}}(E) = H_{M_{\gD_{sf}}} = \text{ the Serre subgroup }
\end{equation}
and also that for each integer $N$ we have
\begin{equation} \label{split&stable}
G_N(E) = \begin{cases} \pi_{N,M_{\gD_{sf}}}^{-1}(H_{M_{\gD_{sf}}}) & \text{ if }
M_{\gD_{sf}} \mid N \\                      
GL_2(\mbz/N\mbz) & \text{ otherwise, }
\end{cases}
\end{equation}
where $\pi_{N,M_{\gD_{sf}}}$ denotes the natural projection
\[
GL_2(\mbz/N\mbz) \longrightarrow GL_2(\mbz/M_{\gD_{sf}}\mbz).
\]
We make use of the following technical lemma.
\begin{lemma} \label{stable} 
Let $N > 1$ be any integer which is
divisible by some prime $p \geq 5$ and by some prime $p < 5$.  Write
\[
N = N_1 \cdot N_2
\]
where $N_1$ is not divisible by any prime $p \geq 5$ and $N_2$ is not divisible
by any prime $p < 5$.
Suppose that $G_a \subset GL_2(\mbz/N\mbz)$ is a subgroup such that
\[
G_a \cap SL_2(\mbz/N\mbz) = (GL_2(\mbz/N\mbz))'.
\]
Finally, assume $G_b \subset G_a$ is a subgroup for which the canonical 
maps
\begin{equation} \label{Gbsurjections}
G_b \twoheadrightarrow GL_2(\mbz/N_1\mbz) \quad \text{ and } \quad G_b \twoheadrightarrow GL_2(\mbz/N_2\mbz)
\end{equation}
as well as the determinant map
\[
\det : G_b \twoheadrightarrow (\mbz/N\mbz)^*
\]
are surjections.  Then $G_b = G_a$.
\end{lemma}
\begin{proof}
Write $N_1  = 2^r 3^s$.  By \eqref{Gbsurjections}, we find by taking commutators that
\[
G_b' \twoheadrightarrow (GL_2(\mbz/N_1\mbz))' = (\ker\ve \cap
SL_2(\mbz/2^r\mbz))\times SL_2(\mbz/3^s\mbz)
\]
and
\[
G_b' \twoheadrightarrow (GL_2(\mbz/N_2\mbz))' = SL_2(\mbz/N_2\mbz).
\]
are also surjections.  We are now in a position to apply the Goursat lemma:
\begin{lemma}
\label{goursat} Let $G_1$ and $G_2$ be groups.  Denote by $\pi_i :
G_1 \times G_2 \rightarrow G_i \; (i = 1,2)$ the projection map.
Suppose that $G \subseteq G_1 \times G_2$ is a subgroup such that
$\pi_i(G) = G_i$ for $i=1,2$ and define
\[
H_1 = \pi_1(G\cap(G_1\times\{e_2\})) \; \text{ and } \; H_2 =
\pi_2(G\cap(\{e_1\}\times G_2)).
\]
Then,
\[
G_1/H_1 \simeq G_2/H_2
\]
and the graph of this isomorphism is induced by $G$.
\end{lemma}

We apply the lemma with $G_1 = (\ker\ve \cap
SL_2(\mbz/2^r\mbz))\times SL_2(\mbz/3^s\mbz)$, $G_2 =
SL_2(\mbz/N_2\mbz)$, and $G = G_b'$ and conclude that $(\ker\ve \cap
SL_2(\mbz/2^r\mbz))\times SL_2(\mbz/3^s\mbz)$ and $SL_2(\mbz/N_2\mbz)$ have a
common quotient group $Q$. If $Q$ is nontrivial then it has a
nontrivial simple quotient $Q_s$.  Since Lemma \ref{jordanhoelder} continues to hold
with $GL_2$ replaced by $SL_2$, we see that $Q_s$ must be abelian.  Since $(SL_2(\mbz/N_2\mbz))' = SL_2(\mbz/N_2\mbz)$, we
conclude that $Q_s = 1$. This shows that $Q$ was trivial to begin
with.  We conclude that $G_1 = H_1$ and $G_2 = H_2$, i.e. that
\[
(GL_2(\mbz/N_1\mbz))'\times\{1\} \subset G_b' \quad \textrm{ and } \quad \{ 1 \}
\times (GL_2(\mbz/N_2\mbz))' \subset G_b',
\]
which implies that 
\[
G_b' = (GL_2(\mbz/N\mbz))'.
\]
But now from the exact sequence
\[
1 \rightarrow (GL_2(\mbz/N\mbz))' \rightarrow G_a \rightarrow (\mbz/N\mbz)^* \rightarrow 1
\]
and
\[
\det : G_b \twoheadrightarrow (\mbz/N\mbz)^*
\]
we conclude that $(GL_2(\mbz/N\mbz))'G_b = G_a$.  So since $(GL_2(\mbz/N\mbz))' \subset G_b$, we have $G_b =
G_a$.
\end{proof}
Now suppose that $E$ is an elliptic curve over $\mbq$ which
satisfies the hypotheses of Lemma \ref{serrecurvelemma}.  First we use
Lemma \ref{stable} to show \eqref{HME}.  This is done by applying
the lemma with $N = M_{\gD_{sf}}$.  Set $G_a = H_{M_{\gD_{sf}}}$ and $G_b = G_{M_{\gD_{sf}}}(E)$.  Write
\[
M_{\gD_{sf}} = 2^r3^s M'
\]
where $M' > 1$ and is co-prime to $6$.  To
see that
\begin{equation} \label{determinantconditiononB} 
H_{M_{\gD_{sf}}} \cap SL_2(\mbz/M_{\gD_{sf}}\mbz) = (\ker
\ve \cap SL_2(\mbz/2^r3^s\mbz)) \times SL_2(\mbz/M'\mbz),
\end{equation}
we argue that the containment ``$\subseteq$'' follows from the
definition of $H_{M_{\gD_{sf}}}$.  To see the reverse containment, we use isomorphism
theorems from group theory and count:
\[
\frac{|H_{M_{\gD_{sf}}}|}{|H_{M_{\gD_{sf}}} \cap SL_2(\mbz/M_{\gD_{sf}}\mbz)|} =
\frac{|GL_2(\mbz/M_{\gD_{sf}}\mbz)|}{|SL_2(\mbz/M_{\gD_{sf}}\mbz)|}.
\]
Now since the index $[GL_2(\mbz/M_{\gD_{sf}}\mbz) : H_{M_{\gD_{sf}}}] = 2$ we see
that 
\[
[SL_2(\mbz/M_{\gD_{sf}}\mbz) : H_{M_{\gD_{sf}}} \cap SL_2(\mbz/M_{\gD_{sf}}\mbz)] = 2.  
\]
Since the index $[SL_2(\mbz/M_{\gD_{sf}}\mbz) : (\ker \ve \cap SL_2(\mbz/2^r3^s\mbz)) \times
SL_2(\mbz/M'\mbz)]$ is obviously equal to 2, we conclude that
\eqref{determinantconditiononB} holds.

We next verify the surjectivity conditions
\begin{equation} \label{surjectivity}
G_{M_{\gD_{sf}}}(E) \twoheadrightarrow GL_2(\mbz/2^r3^s\mbz) \textrm{ and }
G_{M_{\gD_{sf}}}(E) \twoheadrightarrow GL_2(\mbz/M'\mbz).
\end{equation}
If the first map is not surjective, then $E$ has some minimal
exceptional number $d$ which divides $2^r3^s$.  By Lemma
\ref{minimal}, we conclude that $d \in \{2,3,4,8, M_{\gD_{sf}(E)} \}$, where
if $d=8$ there is a specific type of exceptional subgroup
occurring at $8$.  However, the assumptions on $E$ in Lemma
\ref{serrecurvelemma} preclude these possibilities  Therefore
$G_{M_{\gD_{sf}}}(E) \twoheadrightarrow GL_2(\mbz/2^r3^s\mbz)$ is
surjective. Similarly, if $G_{M_{\gD_{sf}}}(E) \rightarrow
GL_2(\mbz/M'\mbz)$ is not surjective, then $E$ has some minimal
exceptional $d$ dividing $M'$.  By Lemma \ref{minimal}, we see that
$d$ must either be a prime number or $d = M_{\gD_{sf}(E)}$.  The assumption
that $E$ has no exceptional primes precludes the first
possibility, and since $2 \nmid M'$, $d$ cannot be even and so $d$
cannot be equal to $M_{\gD_{sf}(E)}$, which is always even.  We have verified
the conditions \eqref{surjectivity}.  Finally, the surjectivity of
\[
\det : G_{M_{\gD_{sf}}} \twoheadrightarrow (\mbz/M_{\gD_{sf}}\mbz)^*
\]
is the surjectivity of the Weil pairing.  By Lemma \ref{stable},
we conclude that
\[
G_{M_{\gD_{sf}}}(E) = H_{M_{\gD_{sf}}}.
\]

Now we verify \eqref{split&stable}.  First let $N$ be any positive
integer and suppose $M_{\gD_{sf}} \nmid N$.  If $G_N(E) \subsetneq
GL_2(\mbz/N\mbz)$ then $E$ has some minimal exceptional number $d$
dividing $N$.  Clearly $d$ cannot be equal to the Serre number
$M_{\gD_{sf}}$, so again by Lemma \ref{minimal} and the assumptions on $E$
from Lemma \ref{serrecurvelemma} we arrive at a contradiction.  Thus, if
$M_{\gD_{sf}} \nmid N$ we have
\[
G_N(E) = GL_2(\mbz/N\mbz).
\]
Now suppose $M_{\gD_{sf}} \mid N$.  We apply Lemma \ref{stable} with $G_a =
\pi_{N,M_{\gD_{sf}}}^{-1}(H_{M_{\gD_{sf}}})$, and $G_b = G_N(E)$.  The verification of
the conditions on $G_a$ and $G_b$ are almost identical to those done
in the previous paragraph, so we omit them.  We conclude in this
case that $G_N(E) = \pi_{N,M_{\gD_{sf}}}^{-1}(H_{M_{\gD_{sf}}})$.

We have shown that for any elliptic curve satisfying the
hypotheses of Lemma \ref{serrecurvelemma} that \eqref{HME} and
\eqref{split&stable} hold, and so our proof of Lemma
\ref{serrecurvelemma} is now complete.

\section{Almost all elliptic curves are Serre curves.} \label{aaecasc}

We now show how Lemma \ref{serrecurvelemma} and Theorem \ref{CHEBOTAREV}
together imply Theorem \ref{main}.

For $N \in \{ 4,6,8,9,12,24 \}$ define
\[
\ve_N(X) := \begin{cases} \{E \in C(X) : E \textrm{ is minimal
exceptional at } N \} & \textrm{ if } N \in \{4,9\} \\
\{E \in C(X) : [GL_2(\mbz/8\mbz) : G_8(E)] = 2  & \textrm{ if } N = 8 \\
\{E \in C(X) : G_N(E) \subseteq H_N \} & \textrm{ if } N \in
\{6,12,24\}. \end{cases}
\]
By Lemma \ref{serrecurvelemma}, we have that the set of non-Serre curves
$C_{nS}(X) \subset C(X)$ satisfies
\[
C_{nS}(X) \subseteq \ve(X) \cup \left( \bigcup_{N \in \{ 4,6,8,9,12,24 \}} \ve_N(X)
\right)
\]
By \eqref{almostall}, to prove Theorem \ref{main} it suffices to estimate the sets
$\ve_N(X)$.
\begin{definition}
Let $W$ be any integer and let $(t,d) \in (\mbz/W\mbz)^2$ be any
pair of integers modulo $W$ with $d \in (\mbz/W\mbz)^*$.  Suppose
that $G \subseteq GL_2(\mbz/W\mbz)$ is any subgroup.  We say that
\emph{$G$ represents the pair $(t,d)$} if there is a matrix $g \in
G$ satisfying
\[
\textrm{tr}\,(g) = t, \quad \det(g) = d.
\]
\end{definition}
The next two lemmas guarantee that when an elliptic curve fails to be a Serre curve by
being exceptional at $N$, there must be some pair $(t,d)$ not represented by $G_N(E)$.
\begin{lemma} \label{MElemma}
Let $W > 2$ be any positive integer and let
\[
\chi : GL_2(\mbz/W\mbz) \rightarrow \{\pm 1\}
\]
be any nontrivial real character.  Suppose that $G \subseteq \ker \chi$ is
any subgroup.  Then there exist integers $t$ and $d$ modulo $W$ so
that the pair $(t,d)$ is not represented by $G$.
\end{lemma}
\begin{proof}
Lemma \ref{commutator} implies that $\chi$ is either of the form $\gd \circ \det$ or $\ve \cdot \gd \circ \det$, 
where $\gd$ is nontrivial and the second possibility may only occur if $W$ is even.  Choose any $f \in (\mbz/W\mbz)^*$ with
\[
\gd(f) = -1
\]
and set $(t,d) = (1,f)$.
\end{proof}
\begin{lemma} \label{4lemma}
Let $p=2$ or $3$ and suppose $G \subseteq GL_2(\mbz/p^2\mbz)$ is a
subgroup which represents every trace-determinant pair $(t,d) \in
(\mbz/p^2\mbz)\times(\mbz/p^2\mbz)^*$ and which surjects onto
$GL_2(\mbz/p\mbz)$.  Then, $G = GL_2(\mbz/p^2\mbz)$.
\end{lemma}
\begin{proof}
We consider the intersection
\[
G \cap K
\]
of $G$ with $K = $ the kernel of the projection
\[
 GL_2(\mbz/p^2\mbz) \longrightarrow GL_2(\mbz/p\mbz).
\]
Our goal is to show that $G$
actually contains $K$.  From here we divide the argument into
cases, according to whether $p$ is $2$ or $3$.

\textbf{Case: $p = 3$.}  Under the given hypothesis, we may find a
matrix $g \in G$ with $\textrm{tr}\,g=3$ and $\det g = 1$.  Such a
matrix must have the form
\[
X + 3Y, \quad X \in \{ \begin{pmatrix} 0 & 1 \\ 2 & 0
\end{pmatrix} , \begin{pmatrix} 0 & 2 \\ 1 & 0
\end{pmatrix} , \begin{pmatrix} 1 & 1 \\ 1 & 2
\end{pmatrix}, \begin{pmatrix} 1 & 2 \\ 2 & 2
\end{pmatrix} , \begin{pmatrix} 2 & 1 \\ 1 & 1
\end{pmatrix}, \begin{pmatrix} 2 & 2 \\ 2 & 1
\end{pmatrix} \},
\]
and the (mod $3$) coefficients of the matrix $Y = \begin{pmatrix}
a & b \\ c & d \end{pmatrix}$ satisfy the conditions
\begin{align*}
a+d=1, \; b-c=1 \textrm{ if } X &= \begin{pmatrix} 0 & 1 \\ 2 & 0
\end{pmatrix} \\
a+d=1, \; b-c=2 \textrm{ if } X &= \begin{pmatrix} 0 & 2 \\ 1 & 0
\end{pmatrix} \\
a+d=0, \; a+b+c-d=0 \textrm{ if } X &= \begin{pmatrix} 1 & 1 \\ 1
& 2 \end{pmatrix} \\
a+d=0, \; a-b-c-d=2 \textrm{ if } X &= \begin{pmatrix} 1 & 2 \\
2 & 2 \end{pmatrix} \\
a+d=0, \; a-b-c-d=0 \textrm{ if } X &= \begin{pmatrix} 2 & 1 \\ 1
& 1 \end{pmatrix} \\
a+d=0, \; a+b+c-d=1 \textrm{ if } X &= \begin{pmatrix} 2 & 2 \\ 2
& 1 \end{pmatrix} .
\end{align*}
In each case, the first equation comes from the trace condition on
$g$ and the second one comes from the determinant condition.  One
computes:
\[
(X+3Y)^4 \equiv I + 3X \mod 9.
\]
Since in this case the discriminant $t^2-4d = 5$ is nonzero modulo
$3$ we see by Lemma \ref{description} that all six of the matrices $X$, when reduced modulo $3$, are
$GL_2(\mbz/3\mbz)$-conjugate to one another.  From this and the
fact that the various $X$ span the $\mbz/3\mbz$-vector space
$M_2(\mbz/3\mbz)$ we conclude that
\[
G \cap K = I + 3M_2(\mbz/3\mbz),
\]
i.e. that $K \subseteq G$, and so $G=GL_2(\mbz/9\mbz)$ in this
case.

\textbf{Case: $p = 2$.}  The proof in this case is similar.  Pick
$g \in G$ with $\text{tr}\,g = 2$ and $\det g = -1$.  Then $g$
must have the form
\[
g = X + 2 Y , \quad X \in \{ \begin{pmatrix} 0 & 1 \\ 1 & 0
\end{pmatrix} , \begin{pmatrix} 1 & 1 \\ 0 & 1
\end{pmatrix}, \begin{pmatrix} 1 & 0 \\ 1 & 1
\end{pmatrix} \}, \quad Y = \begin{pmatrix} a & b \\ c & d
\end{pmatrix}
\]
where the (mod $2$) coefficients of the matrix $Y$ satisfy the
conditions
\begin{align*} a+d = 1, \, b+c = 0 & \text{ if } X = \begin{pmatrix} 0 & 1 \\ 1 & 0
\end{pmatrix} \\
a+d = 0, \, a+c+d = 1& \text{ if } X = \begin{pmatrix} 1 & 1 \\ 0
& 1 \end{pmatrix} \\
a+d = 0, \, a+b+d = 1& \text{ if } X = \begin{pmatrix} 1 & 0 \\ 1
& 1
\end{pmatrix}. 
\end{align*}
(The possibility $X = I + 2Y$ is eliminated since the conditions
on the coefficients of $Y$ in that case read $a+d = 0, \, a+d =
1$.)  One computes:
\[
(X + 2Y)^2 \equiv I + 2X \mod 4.
\]
After conjugating by preimages of elements of $GL_2(\mbz/2\mbz)$,
one concludes that
\[
G \cap K \supseteq \{ I + 2\cdot \, \textrm{span}\, \left\{ \begin{pmatrix} 1 & 1 \\
0 & 1
\end{pmatrix}, \begin{pmatrix} 1 & 0 \\ 1 & 1
\end{pmatrix} \right\} \}.
\]
Playing the same game with $t = 0$ and $d = 1$, one sees that in
fact
\[
G \cap K \supseteq \{ \text{ traceless matrices } \}.
\]
Now $G$ can be at worst an index $2$ subgroup of
$GL_2(\mbz/4\mbz)$.  However, if $G$ is indeed a proper subgroup of index $2$, we
may apply Lemma \ref{MElemma} and arrive at a contradiction.  This concludes the proof in this case.
\end{proof}
Lemmas \ref{MElemma} and \ref{4lemma} imply the
following corollary.
\begin{corollary}
For $N \in \{ 4,6,8,9,12,24 \}$, we have
\[
\ve_N(X) = \bigcup_{(t,d) \in (\mbz/N\mbz)\times(\mbz/N\mbz)^*}
\ve_{N,(t,d)}(X),
\]
where
\[
\ve_{N,(t,d)} := \{ E \in \ve_N(X) : \text{ $(t,d)$ is not
represented by } G_N(E) \}.
\]
\end{corollary}

\begin{lemma} \label{ve23bound}
For each $N \in \{ 4,6,8,9,12,24 \}$, we have
\[
| \ve_N(X) | \ll N^{16} X^6 \max_d \pi(X;N,d)^{-2},
\]
with an absolute implied constant.
\end{lemma}

This lemma and its proof are analogous to Lemma $5$ of \cite{duke}, whose statement contains
a small error: the ``$\ll X^6 \pi(X;N,d)^{-2}$'' should be replaced by ``$\ll N^4 X^6 \max_d \pi(X;N,d)^{-2}$.''

Using the Siegel-Walfisz theorem (c.f. \cite{davenport}, p. 133), which gives $\pi(X;N,d) \gg (\varphi(N))^{-1} \pi(X)$,
Theorem \ref{main} follows.

\section{$N=4$ occurs as a minimal exceptional number.} \label{N=4}
If $N=4$ or $9$, the argument given in section
\ref{minimalexceptional} is invalid since we may not conclude that
\eqref{traceless} holds. In fact, there is a subgroup $H \subset
GL_2(\mbz/4\mbz)$ of index four which satisfies conditions
\eqref{ontod} and \eqref{ontodet}. We now describe $H$ and
demonstrate an infinite family of non-isomorphic elliptic curves
$E$ for which $G_4(E) = H$.  Elkies \cite{elkies} has recently exhibited similar
examples for $N=9$.

First, a geometric description of $H$:  Let $L$ be a complex
lattice and let $L[4]$ denote the 4-torsion of $\mbc/L$.  By
choosing a basis, we may identify $L[4]$ with $(\mbz/4\mbz)^2$.
Let $l_1,l_2 \dots, l_6$ denote the lines through the origin in
$L[4]$.  More precisely, define the equivalence relation on $L[4]$
by declaring $u \sim u'$ exactly if $u' = \gl u$ for some $\gl \in
(\mbz/4\mbz)^* = \{\pm1\}$ and denote the resulting equivalence
classes by $l_1,l_2 \dots, l_6$.  Since the Weierstrass
$\wp$-function is even, the association $l_i = [u] \mapsto \wp
(u)$ identifies $\mbp ((\mbz/4\mbz)^2) := \{l_1,l_2,\dots,l_6 \}$
with $E[4]_x$, the set of $x$-coordinates of the 4-torsion of
$E=E_L$, the elliptic curve associated to the lattice $L$.  This
correspondence identifies the Galois group of $\mbq(E[4]_x)$ over
$\mbq$ with a subgroup of $PGL_2(\mbz/4\mbz)$.

We may extend the natural action of $PGL_2(\mbz/4\mbz)$ on $\mbp
((\mbz/4\mbz)^2)$ to obtain a $PGL_2(\mbz/4\mbz)$ action on the
set
\[
S := \{ \{ \{l_{i_1},l_{i_2},l_{i_3}\} ,
\{l_{i_4},l_{i_5},l_{i_6}\} \} : \text{ all } i_j \in
\{1,2,\dots,6\} \text{ are distinct}\}.
\]
This action is not transitive.  The size 10 set $S$ decomposes
into two orbits $S_1$ and $S_2$ of sizes 4 and 6, respectively. To
describe these sets, one needs to define an ``addition relation''
on $\mbp ((\mbz/4\mbz)^2)$.  If $l_1$, $l_2$, and $l_3$ are lines
in $\mbp ((\mbz/4\mbz)^2)$, we say that
\[
l_1 + l_2 = l_3
\]
exactly when, for some choice of representatives $u_i \in l_i$ we
have
\[
u_1 + u_2 = u_3.
\]
(Note:  This is a \emph{relation}, not a well-defined operation.
For example, $[(1,0)] + [(0,1)] = [(1,1)]$ and $[(1,-1)]$.) Then
the two orbits are defined by
\[
S_1 := \{ \{ \{l_{i_1},l_{i_2},l_{i_3}\} ,
\{l_{i_4},l_{i_5},l_{i_6}\} \} \in S : l_{i_1} + l_{i_2} = l_{i_3}
\}
\]
and $S_2 = S - S_1$. Fixing any element $r \in S_1$ we define $H_x
= H_x(r) \subset PGL_2(\mbz/4\mbz)$ to be the stabilizer of $r$.
Finally, we define $H = H(r)$ to be the preimage of $H_x$ under
the natural projection $GL_2(\mbz/4\mbz) \rightarrow
PGL_2(\mbz/4\mbz)$.

To find elliptic curves $E$ with $G_4(E) = H$, we reason as
follows:  let $x_1$, $x_2$, \dots, $x_6$ be the elements of
$E[4]_x$.  If $E$ is given in the form
\[
E: \quad y^2 = 4x^3-g_2x-g_3
\]
then the minimal polynomial for $x_1$, $x_2$, \dots, $x_6$ is
given by
\[
f_E(t) = t^6-\frac{5g_2}{4}t^4-5g_3t^3-\frac{5g_2^2}{16}t^2
-\frac{g_2g_3}{4}t+\frac{g_2^3-32g_3^2}{64}.
\]
The set $S_1$ defined above corresponds to the set of numbers
\[
X_1 := \{ (x_{i_1} + x_{i_2} + x_{i_3})(x_{i_4} + x_{i_5} +
x_{i_6}) : (x_{i_1},y_{i_1}) \oplus (x_{i_2},y_{i_2}) =
(x_{i_3},y_{i_3}) \},
\]
where $\oplus$ refers to the addition law on $E$.  $X_1$ is a set
of four complex numbers which satisfy the (generically
irreducible) polynomial
\[
f_{1,E}(t) = t^4+3g_2t^3+\frac{27g_2^2}{8}t^2+
(\frac{-37g_2^3}{16}+108g_3^2)t+\frac{81g_2^4}{256}.
\]
We note that $G_4(E) \subseteq$ some $H(r)$ whenever $f_{1,E}(t)$
has a linear factor over $\mbq$. Let $s \in \mbq$ and denote by
$E_s$ the elliptic curve given by the equation
\[
y^2 :=
4x^3+\frac{16s^2+56s+81}{3s}x+\frac{(16s^2+56s+81)^2(-1+4s)}{864s^2}.
\]
It may be checked that for each $s$, $f_{1,E_s}(t)$ is divisible
by $t+27+\frac{56}{3}s+\frac{16}{3}s^2$ and that
\[
\text{Gal}(\mbq(E_s[4])/\mbq(s)) \simeq H.
\]
The discriminant is computed to be
\[
\gD(E_s) = -\frac{(16s^2+56s+81)^3(4s+3)^4}{27648s^4},
\]
and the $j$-invariant is
\[
j(E_s) = \frac{1769472s}{(4s+3)^4}.
\]
In particular, if we apply the Hilbert irreducibility criterion, we see that there
are infinitely many non-isomorphic curves $E_s$, each with Galois
group $G_4 \simeq H$.

\vspace{2in}

\begin{center}
Centre de Recherches Math\'{e}matiques \\
Universit\'{e} de Montr\'{e}al \\
P.O. Box 6128, \\
Centre-ville Station \\
Montr\'{e}al, Qu\'{e}bec  H3C 3J7, Canada. \\
E-mail:  jones@dms.umontreal.ca
\end{center}

\end{document}